\documentclass[letterpaper,11pt]{amsart}

\usepackage{amssymb}
\usepackage{amscd}
\usepackage{eucal}
\usepackage{epsfig}
\usepackage[matrix,arrow,curve,tips]{xy}
\CompileMatrices

\newcommand{\includefigure}[3]{%
\begin{figure}[#3]
\begin{center}
\epsfig{file=#2} \\ \caption{\label{fig:#1}}
\end{center}
\end{figure}}

{\begin{list}{{\rm(\alph{enumi}) }}{\usecounter{enumi}

\leftmargin0cm \labelsep0cm \rightmargin0cm \parsep1em \listparindent0em
\itemsep-1em \topsep1em \parskip1em \setlength{\labelwidth}{\fill}}}
{\parskip0em \end{list}}

\newcommand{\R}{\mathbb{R}}
\newcommand{\Z}{\mathbb{Z}}

\newcommand{\C}{\mathbb{C}}

\newcommand{\iso}{\cong}           
\newcommand{\htp}{\simeq}          

\newcommand{\leftsc}{\langle}
\newcommand{\rightsc}{\rangle}

\newcommand{\id}{\mathrm{id}}

\newcommand{\im}{\mathrm{im}}
\renewcommand{\ker}{\mathrm{ker}}
\newcommand{\coker}{\mathrm{coker}}

\newcommand{\Hom}{\mathrm{Hom}}

\renewcommand{\o}{\omega}

\newcommand{\Aut}{\mathrm{Aut}}

\newcommand{\Symp}{\mathrm{Symp}}
\newcommand{\Diff}{\mathrm{Diff}}
\newcommand{\Lag}{\mathrm{Lag}}
\SelectTips{cm}{}

\renewcommand{\subsection}{\vspace{1ex}
\hspace{-\parindent}%
\stepcounter{subsection}%
({\sc \thesection\alph{subsection}})}
\theoremstyle{plain}

\newtheorem{thm}{Theorem}[section]
\newtheorem{theorem}[thm]{Theorem}

\newtheorem{corollary}[thm]{Corollary}
\newtheorem{lemma}[thm]{Lemma}

\newtheorem{prop}[thm]{Proposition}

\newtheorem*{conjecture*}{Conjecture}

\newtheorem*{question*}{Question}

\newtheorem{definition}[thm]{Definition}

\newtheorem*{definitions*}{Definitions}

\newtheorem*{rem*}{Remark}
\newtheorem{remark}[thm]{Remark}
\newtheorem*{remark*}{Remark}

\newtheorem*{remarks*}{Remarks}
\newtheorem*{example*}{Example}
\newtheorem{example}[thm]{Example}
\newtheorem*{examples*}{Examples}

\newtheorem*{notation*}{Notation}
\newtheorem*{convention*}{Convention}
\newtheorem*{conventions*}{Conventions}
\newtheorem*{note*}{Note}
\newtheorem*{exercise*}{Exercise}

\theoremstyle{remark}
\newtheorem*{references}{References}
\newtheorem*{acknowledgments}{Acknowledgments}

\newenvironment{myitemize}%
{\begin{itemize}\itemsep0.5em} {\end{itemize}}
\newcommand{\Sympe}{\Symp^e}
\newcommand{\crit}{\mathrm{crit}}

\newcommand{\Phirel}{\Phi_{rel}}

\newcommand{\A}{\mathcal A}
\newcommand{\B}{\mathcal B}
\newcommand{\Aoplus}{{\A^{\oplus}}}
\newcommand{\Cone}{\mathrm{Cone}}

\newcommand{\tL}{\widetilde{L}}
\newcommand{\tphi}{\tilde{\phi}}
\newcommand{\ttau}{\tilde{\tau}}

\newcommand{\J}{\mathbf{J}}

\renewcommand{\Lag}{{\mathit{Lag}}^\rightarrow}

\renewcommand{\hom}{\mathit{hom}}
\newcommand{\Tw}{\mathrm{Tw}}
\newcommand{\Ob}{\mathrm{Ob}}
\newcommand{\JJ}{\mathcal{J}}
\newcommand{\RR}{\mathcal{R}}
\renewcommand{\SS}{\mathcal{S}}
\newcommand{\CC}{\mathcal{C}}
\title{Vanishing cycles and mutation}
\author{Paul Seidel}
\date{July 16, 2000.}

\begin{document}
\begin{abstract}
Using Floer cohomology, we establish a connection between
Picard-Lefschetz theory and the notion of mutation of exceptional
collections in homological algebra.
\end{abstract}
\maketitle

\section{Introduction}

This talk is about symplectic aspects of Picard-Lefschetz theory, and the role
of Floer cohomology in that context. I have relied on two main sources for
inspiration, which are the ideas of Donaldson on vanishing cycles
\cite{donaldson98} respectively those of Kontsevich, partly in collaboration
with Barannikov, on mirror symmetry for Fano varieties \cite{kontsevich98}. On
a more technical level, the basis is provided by Fukaya's work on Floer
cohomology \cite{fukaya93}. I have tried to be as concrete as possible: not
only are all objects mentioned rigorously defined, but they can be explicitly
computed, and the assertions made about them checked, in many examples. This
hands-on approach has its drawbacks, one of which will be mentioned after the
summary of contents.

The basic geometric notion is that of an exact Morse fibration, a symplectic
analogue of a holomorphic Morse function. One way of analyzing such fibrations
is through the vanishing cycles in a fibre. When the base is a disc, one
conveniently makes a choice of a finite family of vanishing cycles, forming a
so-called distinguished basis. Any two such bases are connected by a sequence
of Hurwitz moves. These are familiar concepts except that here vanishing cycles
are considered as Lagrangian submanifolds, rather than only as homology
classes.

Apart from their role as geometric objects worthy of study on their own, exact
Morse fibrations are also relevant to Floer theory since, when equipped with
suitable Lagrangian boundary conditions, they provide homomorphisms between
Floer cohomology groups. As an application of these new maps we construct a
long exact sequence analogous to that of Floer in gauge theory.

After that we return to exact Morse fibrations over a disc. The goal is to
associate to each such fibration a triangulated category. An additional
assumption is necessary: for simplicity, let's say that the total space of the
fibration has zero first Chern class and zero first Betti number. The
construction proceeds in several steps: first one chooses a distinguished basis
of vanishing cycles. From that one obtains a Fukaya-type $A_\infty$-category,
unique up to quasi-isomorphism; our invariant is the derived category of this.
It resembles derived categories of coherent sheaves on some Fano varieties, in
that it is generated by an exceptional collection. The main point is to show
that Hurwitz moves of the distinguished basis correspond to mutations of the
exceptional collection, which leave the derived category unchanged. The
explicit nature of mutations also allows them to be used for concrete
computations of Floer cohomology.

We have to admit that our triangulated categories are only well-defined in a
weak sense: different choices made during the construction lead to equivalent
categories, but it has not been proved that the equivalences are canonical up
to isomorphism (to have a completely satisfactory theory, it would further be
necessary to establish coherence relations between functor isomorphisms).
Rather than trying to do these improvements, it seems better to look for an
alternative definition bypassing the choice of distinguished basis. The
approach envisaged by Kontsevich is of this kind, but more work would be needed
to put it on a rigorous footing.

\begin{acknowledgments}
Thanks go to Ivan Smith for reading the manuscript and suggesting several
improvements.
\end{acknowledgments}

\section{Picard-Lefschetz theory\label{sec:picard-lefschetz}}

\subsection{}
Let $(M,\o,\theta)$ be an exact symplectic manifold of dimension $2n$. This
means that $M$ is compact with boundary, $\o \in \Omega^2(M)$ is a
symplectic form, and $\theta \in \Omega^1(M)$ satisfies $d\theta = \o$. We will
consider only symplectic automorphisms $\phi$ of $M$ which are equal to the
identity near $\partial M$. Such a $\phi$ is called exact if
$[\phi^*\theta-\theta] \in H^1(M,\partial M;\R)$ is zero. The exact symplectic
automorphisms form a subgroup $\Sympe(M) \subset \Symp(M)$. Note that any
isotopy within this subgroup is Hamiltonian.

Let $S$ be a smooth manifold with boundary. An exact symplectic fibration over
$S$ consists of data $(E,\pi,\Omega,\Theta)$ as follows. $\pi: E \rightarrow S$
is a proper differentiable fibre bundle whose fibres are $2n$-dimensional
manifolds with boundary. This means that $E$ itself is a manifold with
codimension two corners, with $\partial E = \partial_vE \cup \partial_hE$
consisting of two faces: $\partial_vE = \pi^{-1}(\partial S)$, while
$\pi|\partial_hE:\partial_hE \rightarrow S$ is again a differentiable fibre
bundle. $\Omega \in \Omega^2(E)$ is closed, its vertical part
$\Omega|\ker(D\pi)$ nondegenerate at every point, and $\Theta \in \Omega^1(E)$
satisfies $d\Theta = \Omega$ . We impose a final condition of triviality near
$\partial_hE$. This means that there should be a neighbourhood $W \subset E$ of
$\partial_hE$ and a diffeomorphism, for some $z \in S$, $\Xi: S \times (W \cap
E_z) \rightarrow W$ lying over $S$, such that $\Xi^*\Omega =
pr_2^*(\Omega|E_z)$ and $\Xi^*\Theta = pr_2^*(\Theta|E_z)$; here $pr_2$ is
projection from $S \times (W \cap E_z)$ to the second factor. Clearly each
fibre $(E_z,\o_z=\Omega|E_z,\theta_z = \Theta|E_z)$ is an exact symplectic
manifold. The form $\Omega$ defines a canonical connection on $\pi: E
\rightarrow S$, with structure group $\Sympe(E_z)$. In fact there is a
bijective correspondence between fibre bundles with structure group
$\Sympe(E_z)$ in the usual sense, and cobordism classes of exact symplectic
fibrations. We denote the parallel transport maps of the canonical connection
by $\rho_c: E_{c(a)} \rightarrow E_{c(b)}$, for $c: [a;b] \rightarrow S$.

\subsection{}
From now on assume that $S$ is two-dimensional and oriented. An {\bf exact
Morse fibration} (this is shorthand for ``exact symplectic fibration with
Morse-type critical points'') over $S$ consists of data
$(E,\pi,\Omega,\Theta,J_0,j_0)$. The properties of $E,\pi,\Omega,\Theta$ are as
before, except that $\pi$ is allowed to have finitely many critical points.
Each fibre may contain at most one of these points, and there should be none at
all on $\partial E$. $J_0$ is an integrable complex structure defined in a
neighbourhood of the set $E^\crit \subset E$ of critical points, and $\Omega$
must be a K{\"a}hler form for it. Similarly $j_0$ is a positively oriented
complex structure on a neighbourhood of the set $S^\crit \subset S$ of critical
values. They should be such that $\pi$ is $(J_0,j_0)$-holomorphic, with
nondegenerate second derivative at each critical point. We will usually denote
exact Morse fibrations by $(E,\pi)$ only.

One thing that needs explaining is why these are supposed to be analogues of
holomorphic functions. For this one considers pairs $(j,J)$ consisting of a
positively oriented complex structure $j$ on $S$ and an almost complex
structure $J$ on $E$, such that $j = j_0$ near $S^\crit$, $J = J_0$ near
$E^\crit$, $\pi$ is $(J,j)$-holomorphic, and $\Omega(\cdot,J\cdot)|\ker(D\pi)$
is symmetric and positive definite everywhere. In this situation we say that
$J$ is compatible relative to $j$. The space of such pairs $(j,J)$ is always
contractible, and in particular nonempty. Moreover, for a fixed pair, by adding
a positive two-form from $S$ one can modify $\Omega$ such that it becomes
symplectic and tames $J$.

Restricting any exact Morse fibration to $S \setminus S^\crit$ yields an exact
symplectic fibration. Before bringing the singular fibres into the picture, we
need some more definitions. Let $(M,\o,\theta)$ be an exact symplectic
manifold. A Lagrangian submanifold $L \subset M$, always assumed to be disjoint
from $\partial M$, is called exact if $[\theta|L] \in H^1(L;\R)$ is zero. A
{\bf framed Lagrangian sphere} is a Lagrangian submanifold $L$ together with an
equivalence class $[f]$ of diffeomorphisms $f: S^n \rightarrow L$. Here
$f_1,f_2$ are equivalent if $f_2^{-1}f_1$ is isotopic to some element of
$O(n+1) \subset \Diff(S^n)$. One can associate to any $(L,[f])$ a Dehn twist
$\tau_{(L,[f])} \in \Symp(M)$ which is unique up to Hamiltonian isotopy. If $L$
is exact, so is the Dehn twist along it. In future, we will often omit the
framing $[f]$ from the notation.

To return to our discussion, let $(E,\pi)$ be an exact Morse fibration. Take a
path $c: [0;1] \rightarrow S$ with $c^{-1}(S^\crit) = \{1\}$ and $c'(1) \neq
0$. Let $x$ be the unique critical point in $E_{c(1)}$. Then the stable
manifold
\[
 B = \{ y \in E_{c(s)}, \, 0 \leq s < 1,
 \text{ with } \lim_{t \rightarrow 1} \,
 \rho_{c|[s,t]}(y) = x\} \cup \{x\}
\]
is a smoothly embedded $(n+1)$-dimensional ball on which $\Omega$ vanishes
identically, and therefore $V = \partial B = B \cap E_{c(0)}$ is an exact
Lagrangian submanifold of $E_{c(0)}$ diffeomorphic to $S^n$. Moreover, $V$ has
a canonical structure of a framed Lagrangian sphere, constructed by first
carrying $V$ by parallel transport to $B \cap E_{c(s)}$, for some $s$ close to
$1$, then projecting orthogonally in local K{\"a}hler coordinates to ${TB}_x$,
and finally projecting radially to the unit sphere in that tangent space. The
composition of these maps is a diffeomorphism $f^{-1}: V \rightarrow S^n$ whose
inverse is the framing. $(V,[f])$ is called the {\bf vanishing cycle}
associated to $c$. A symplectic version of the Picard-Lefschetz theorem says
that for $l$,$c$ as in Figure \ref{fig:doubling}, the monodromy $\rho_l \in
\Sympe(E_{c(0)})$ is isotopic to $\tau_{(V,[f])}$.
\includefigure{doubling}{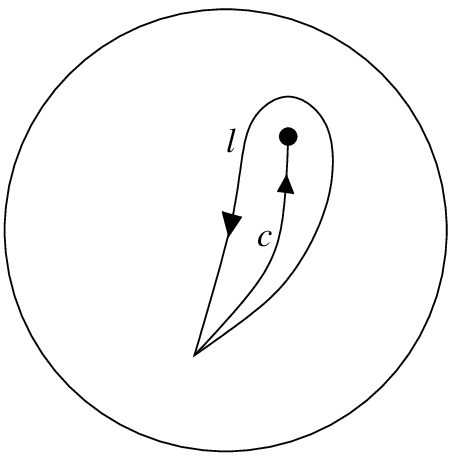}{ht}%

\subsection{}
Suppose now that $S = D$ is the closed unit disc in $\C$. Let $(E,\pi)$ be an
exact Morse fibration over it with $m$ critical values. Let $M$ be the fibre at
some base point $z_0 \in \partial D$, say $z_0 = -i$. An admissible choice of
paths is a family $(c_1,\dots,c_m)$ looking as in Figure \ref{fig:basis},
ordered by their tangent directions at $z_0$. The collection of exact framed
Lagrangian spheres in $M$ arising from them is called a {\bf distinguished
basis} of vanishing cycles. Modifying the paths affects the distinguished basis
in a way which can be determined using the Picard-Lefschetz theorem. The
outcome is encoded in the following abstract notion:
\includefigure{basis}{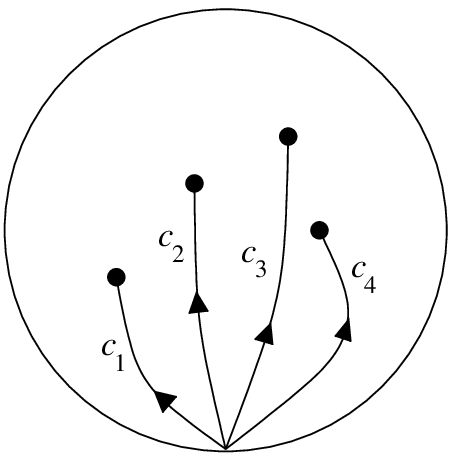}{hb}%

\begin{definition} \label{def:hurwitz}
A Lagrangian configuration in an exact symplectic manifold $M$ is an ordered
family $\Gamma = (L_1,\dots,L_m)$ of exact framed Lagrangian spheres. Two
configurations are Hurwitz equivalent if they can be connected by a sequence of
the following moves and their inverses:
\begin{itemize}
\item
$\Gamma = (L_1,\dots,L_m) \rightsquigarrow
(L_1,\dots,L_{i-1},\phi(L_i),L_{i+1},\dots,L_m)$ for some $1 \leq i \leq m$ and
some $\phi \in \Sympe(M)$ isotopic to the identity;
\item
$\Gamma \rightsquigarrow c\Gamma = (\tau_{L_1}(L_2),\tau_{L_1}(L_3),\dots,
\tau_{L_1}(L_m),L_1)$;
\item
$\Gamma \rightsquigarrow r\Gamma =
(L_1,\dots,L_{m-2},\tau_{L_{m-1}}(L_m),L_{m-1})$.
\end{itemize}
\end{definition}

Thus, the Hurwitz equivalence class of a distinguished basis is an invariant of
the exact Morse fibration $(E,\pi)$. Conversely, from $M$ and that Hurwitz
equivalence class one can reconstruct the fibration up to a suitable notion of
deformation equivalence.

\section{Floer cohomology\label{sec:floer}}

\subsection{}
From now on, any exact symplectic manifold $(M,\o,\theta)$ is assumed to have
contact type boundary, with $\theta|\partial M$ being the contact one-form (the
same condition will be imposed on $\Theta|\partial E_z$, for $E_z$ any fibre of
a Morse fibration). Then for any pair $(L_1,L_2)$ of exact Lagrangian
submanifolds in $M$ there is a well-defined Floer cohomology group
$HF(L_1,L_2)$, which is a finite-dimensional vector space over the field
$\Z/2$. We remind the reader that this is invariant under isotopies of $L_1$ or
$L_2$, satisfies $HF(\phi L_1, \phi L_2) \iso HF(L_1,L_2)$ for any $\phi \in
\Sympe(M)$, and that there is a natural Poincar{\'e} duality $HF(L_1,L_2) \iso
HF(L_2,L_1)^\vee$.

As a warm-up exercise, suppose that we have a compact oriented surface $S$ with
boundary, and an exact Morse fibration $(E^{2n+2},\pi)$ over it. A {\bf
Lagrangian boundary condition} for $E$ is a closed submanifold $Q^{n+1} \subset
\partial_vE \setminus \partial_hE$ such that $\pi|Q: Q \rightarrow \partial S$
is a smooth fibration, satisfying $\Omega|Q = 0$ and $[\Theta|Q] = 0 \in
H^1(Q;\R)$. Then the intersection $Q_z = Q \cap E_z$, for any $z \in \partial
S$, is an exact Lagrangian submanifold in $E_z$, and parallel transport along
$\partial S$ takes these Lagrangian submanifolds into each other. Choose a
complex structure $j$ on $S$ and an almost complex structure $J$ on $E$ which
is compatible relative to $j$, as defined in the previous section. There is a
Gromov type invariant $\Phi(E,\pi,Q) \in \Z/2$ which counts, in the familiar
sense, the number of $(j,J)$-holomorphic sections $u: S \rightarrow E$ with
$u(\partial S) \subset Q$. The exactness assumptions imply that there can be no
bubbles ($J$-holomorphic spheres in a fibre $E_z$, or $J$-holomorphic discs in
$E_z$ with boundary in $Q_z$), so that the definition of the invariant is
technically quite simple.

\begin{example} \label{ex:model-fibration}
Let $L$ be an exact framed Lagrangian sphere in an exact symplectic manifold
$M$. Starting from a standard local model, one can construct an exact Morse
fibration $(E,\pi)$ over $D$ with $E_{z_0} = M$ for $z_0 = -i \in \partial D$,
having exactly one critical point, such that the monodromy around $\partial D$
is $\tau_L$. Because $\tau_L(L) = L$, there is a unique Lagrangian boundary
condition $Q \subset E$ with $Q_{z_0} = L$. $\Phi(E,\pi,Q)$ vanishes because
the expected dimension of the space of $(j,J)$-holomorphic sections is always
odd, hence never zero.
\end{example}

Now let $S$ be as before but with a finite set of marked points $\Sigma \subset
\partial S$. Suppose moreover that around each $\zeta \in \Sigma$ we have
preferred local coordinates, given by an oriented embedding $\psi_\zeta: D^+
\rightarrow S$ of the half-disc $D^+ = D \cap \{\im(z) \geq 0\}$ with
$\psi_\zeta(0) = \zeta$. Let $(E,\pi)$ be an exact Morse fibration over $S^* =
S \setminus \Sigma$ which is trivial near the marked points. This means that we
have a fixed exact symplectic manifold $M$ and preferred embeddings
$\Psi_\zeta: (D^+ \setminus \{0\}) \times M \rightarrow E$ lying over
$\psi_\zeta$, satisfying some obvious conditions concerning $\Omega$,$\Theta$
that we do not care to write down. If $Q \subset E$ is a Lagrangian boundary
condition, there is for each $\zeta \in \Sigma$ a unique pair $L_{\zeta,\pm}$
of exact Lagrangian submanifolds of $M$ such that $\Psi_\zeta^{-1}(Q) = [-1;0)
\times L_{\zeta,-} \, \cup \, (0;1] \times L_{\zeta,+}$. After choosing a
complex structure $j$ on $S^*$ such that the $\psi_\zeta$ become holomorphic,
and an almost complex structure $J$ on $E$ which is compatible relative to $j$
and satisfies some additional conditions with regard to $\Psi_\zeta$, one can
count pseudo-holomorphic sections with suitable behaviour near the marked
points. The outcome is a relative invariant
\begin{equation} \label{eq:phirel}
\Phirel(E,\pi,Q) \in \bigotimes_{\zeta \in \Sigma} HF(L_{\zeta,+},L_{\zeta,-}).
\end{equation}
These invariants satisfy the standard gluing law for a topological quantum
field theory, which one can formulate in two parts as follows. First, if $S$ is
not connected then the relative invariant decomposes into the tensor product of
relative invariants associated to its connected components. Second, suppose
that there are two marked points $\zeta,\zeta' \in \Sigma$ with
\[
L_{\zeta,\pm} = L_{\zeta',\mp}.
\]
One can define a new surface $\overline{S}$ by removing small half-discs around
$\zeta$,$\zeta'$ and gluing together the resulting half-circles, as in Figure
\ref{fig:gluing}. There is a natural set of marked points $\overline{\Sigma}
\subset \partial \overline{S}$ which is inherited from $\Sigma \setminus
\{\zeta,\zeta'\}$. A similar process applied to $(E,\pi)$ constructs a new
exact Morse fibration over $\overline{S} \setminus \overline{\Sigma}$ with
Lagrangian boundary conditions. The gluing rule says that on the level of the
invariants $\Phirel$ this translates into contracting $HF(L_{\zeta,+},
L_{\zeta,-}) \otimes HF(L_{\zeta',+},L_{\zeta',-})$ by Poincar{\'e} duality.
\includefigure{gluing}{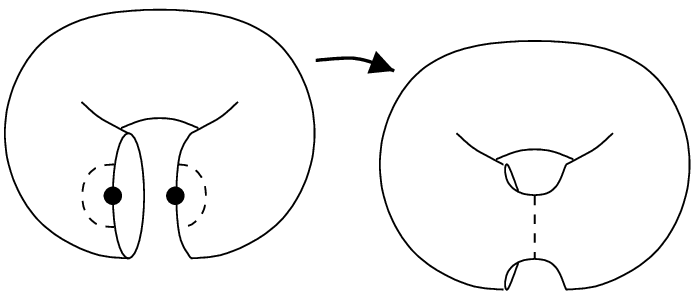}{ht}%

\begin{remark}
The reader is hereby warned of two possible misunderstandings: the gluing
process does not take place along the boundaries $\partial E_z$ of the fibres,
and neither do we glue together two boundary circles of $S$.
\end{remark}

\subsection{}
We next give three examples. Throughout, Poincar{\'e} duality will be used
freely to write the invariants $\Phirel$ in various equivalent ways, e.g.\ as
multilinear maps between Floer cohomology groups.

\begin{myitemize}
\item
Let $M$ be an exact symplectic manifold and $L_1,L_2,L_3 \subset M$ exact
Lagrangian submanifolds. Take $S = D$ and $\Sigma = \{\text{\em three points}\}
\subset \partial S$. Let $I_1,I_2,I_3$ be the three components of $\partial
S^*$, ordered in positive sense. The trivial fibration $E = S^* \times M$ with
Lagrangian boundary conditions $Q = \bigcup_{\nu=1}^3 I_\nu \times L_\nu$
yields a relative invariant which can be written as a map $HF(L_2,L_3) \otimes
HF(L_1,L_2) \rightarrow HF(L_1,L_3)$. This is just a reformulation of the usual
``pair-of-pants'' product. Reversing the orientation of $S$ yields another
relative invariant, which is the Poincar{\'e} dual coproduct $HF(L_1,L_3)
\rightarrow HF(L_2,L_3) \otimes HF(L_1,L_2)$.

\item
Let $M$ be an exact symplectic manifold, $L_1,L_2 \subset M$ exact Lagrangian
submanifolds, and $\phi \in \Sympe(M)$ an automorphism which is isotopic to the
identity. Imagine $S = D$ as being constructed out of a smaller disc $D_0$ and
two other pieces $D_1,D_2 = [0;1]^2$. The marked points are $\Sigma =
\{\zeta_1,\zeta_2\}$, where $\zeta_1$ is obtained by identifying $(0,1) \in
D_1$ with $(0,0) \in D_2$, and $\zeta_2$ is similarly $(1,1) \in D_1$ or $(1,0)
\in D_2$. One can construct an exact symplectic fibration $E_0$ over $D_0$,
which is topologically $D_0 \times M$ but has nontrivial forms $\Omega$ and
$\Theta$, such that the monodromy around $\partial D_0$ is $\phi$. Assemble
$E_0$ and the two trivial fibrations $E_1 = (D_1 \setminus \{(0,1),(1,1)\})
\times M$, $E_2 = (D_2 \setminus \{(0,0),(0,1)\}) \times M$ following the
instructions in Figure \ref{fig:continuation}. This gives an exact Morse
fibration over $S^*$; we equip it with the Lagrangian boundary condition which
is the union of $(0;1) \times \{1\} \times L_1 \subset E_1$ and $(0;1) \times
\{0\} \times L_2 \subset E_2$. The resulting relative invariant is a map
$HF(L_1,L_2) \rightarrow HF(L_1,\phi(L_2))$, in fact just the familiar
``continuation'' isomorphism.
\includefigure{continuation}{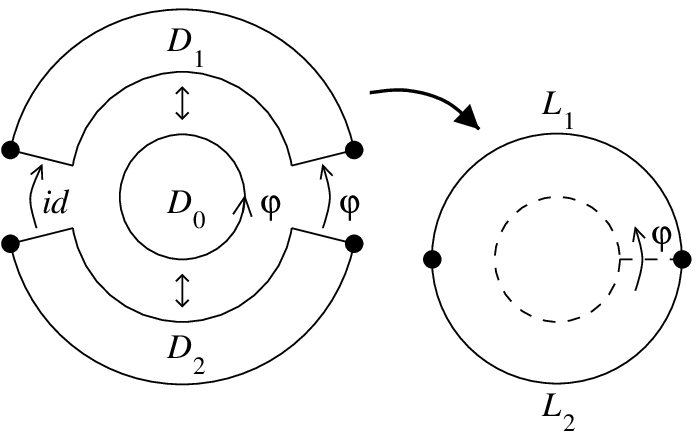}{ht}%

\item
Let $M,L_1,L_2$ be as before, and $L \subset M$ an exact framed Lagrangian
sphere. The central piece $E_0$ in the previous construction can be replaced by
the Morse fibration from Example \ref{ex:model-fibration}. This leads to a map
$HF(L_1,L_2) \rightarrow HF(L_1,\tau_L(L_2))$, which is not an isomorphism in
general. In fact there seems to be no way at all of getting from our formalism
a map in the inverse direction; the reason is the presence of critical points,
which prevents one from reversing the orientation of the base.
\end{myitemize}

\begin{theorem} \label{th:exact-sequence}
Let $(M,\o,\theta)$ be an exact symplectic manifold, $L_1,L_2 \subset M$ exact
Lagrangian submanifolds, and $L \subset M$ an exact framed Lagrangian sphere.
Suppose that $2c_1(M,L) \in H^2(M,L)$ is zero. Then there is a long exact
sequence, with $a$ the pair-of-pants product and $b$ the map defined in the
last example above,
\begin{equation} \label{eq:exact-sequence}
\xymatrix{
 {HF(L,L_2) \otimes HF(L_1,L)} \ar[r]^-{a} &
 {HF(L_1,L_2)} \ar[d]^{b} \\
 & \ar@/^1pc/[ul] {HF(L_1,\tau_L(L_2)).}
 }
\end{equation}
\end{theorem}

To understand how this works, one needs to look at the Floer cochain complexes
$CF$. To simplify, we suppress the dependence of these complexes on various
additional choices, and write $a,b$ for the maps between them inducing the
Floer cohomology maps mentioned above. The central object in the proof is a map
of complexes
\begin{equation} \label{eq:quasi-iso}
\Cone(a) \stackrel{(h,b)}{\longrightarrow} CF(L_1,\tau_L(L_2)).
\end{equation}
Here $h: CF(L,L_2) \otimes CF(L_1,L) \rightarrow CF(L_1,\tau_L(L_2))$ is a
chain homotopy $b \circ a \htp 0$ defined as follows. Consider the exact Morse
fibration with $S = D$, $\Sigma = \{\text{\em three points}\}$, which is
represented schematically, together with its Lagrangian boundary condition, in
Figure \ref{fig:homotopy}. Moving the two leftmost marked points simultaneously
along $\partial D$, in a way that preserves the symmetry of the picture with
respect to the $x$-axis, yields a one-parameter family of fibrations, and $h$
arises from the corresponding parametrized spaces of pseudo-holomorphic
sections. Once \eqref{eq:quasi-iso} has been defined, an argument using the
natural filtration of the Floer complexes by the action functional shows that
it is a quasi-isomorphism; the long exact sequence is an immediate consequence.
\includefigure{homotopy}{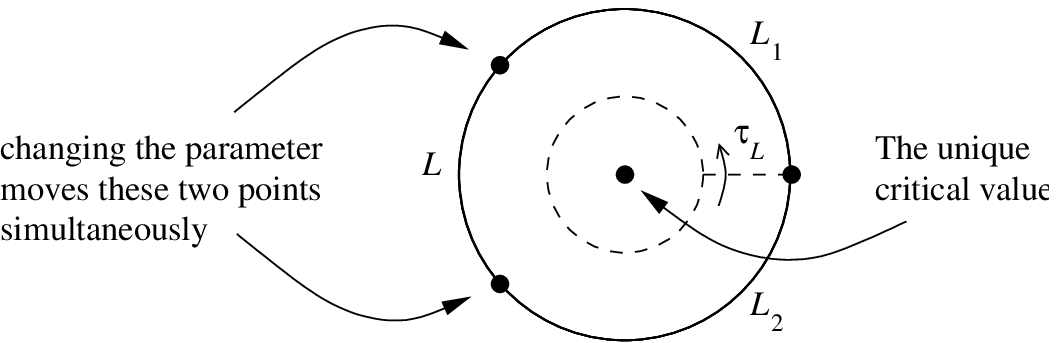}{ht}%

As one can see from this sketch of the argument, the definition of the third
arrow in \eqref{eq:exact-sequence} uses the inverse of the quasi-isomorphism
\eqref{eq:quasi-iso}. Poincar{\'e} duality yields a symmetry of the exact
sequence, which suggests a more direct description of that arrow. Namely, it
should be the composition
\begin{multline*}
c: HF(L_1,\tau_L(L_2)) \rightarrow HF(L,\tau_L(L_2)) \otimes HF(L_1,L) \iso \\
\iso HF(L,L_2) \otimes HF(L_1,L)
\end{multline*}
of the coproduct and the isomorphism $HF(L,\tau_L(L_2)) \iso
HF(\tau_L^{-1}(L),L_2) = HF(L,L_2)$. One can show that this is indeed the same
map as that obtained by inverting \eqref{eq:quasi-iso}; the proof uses
invariants of the same kind as those defining $h$, for two-parameter families
of exact Morse fibrations.

\begin{remark} \label{rem:technical}
The assumption $2c_1(M,L) = 0$ is a technical one. It implies that the space of
$(j,J)$-holomorphic sections in Example \ref{ex:model-fibration} has expected
dimension $2n-1$. From this it follows (except for $n = 1$, which requires a
separate treatment) that generically there is a positive-dimensional space of
these sections $u$ such that $u(z_0)$ is a specific point in $L$. That enters
into the description of the limiting behaviour of sections of the family in
Figure \ref{fig:homotopy} when both moveable points of $\Sigma$ go towards the
fixed one, and through it into the proof that $h$ is a homotopy $b a \htp 0$.
Still, it may be possible to substitute some other argument at this point, and
thereby remove the assumption from Theorem \ref{th:exact-sequence}.
\end{remark}
\includefigure{skein}{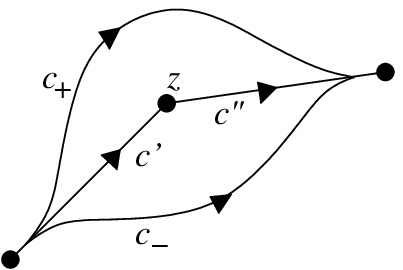}{hb}%

There is an important special case in which Theorem \ref{th:exact-sequence} has
a graphical interpretation, namely when all Lagrangian submanifolds involved
are vanishing cycles. Take an exact Morse fibration $(E,\pi)$ with arbitrary
base $S$. To any path $c: [0;1]
\rightarrow S$ with $c^{-1}(S^\crit) = \{0;1\}$ and $c'(0),c'(1) \neq 0$ one
can associate the Floer cohomology group $HF(V_{c,0},V_{c,1})$ where $V_{c,0},
V_{c,1} \subset E_{c(t)}$, for some $0<t<1$, are the vanishing cycles
associated to $c|[0;t]$ resp.\ $c|[t;1]$. This group is independent of $t$, so
we may write temporarily $HF(c)$ for it. Consider four paths as in Figure
\ref{fig:skein}. As three of them run together for small $t$, we may assume
that their vanishing cycles all lie in the same fibre, called $M$. Then
$V_{c_-,0} = V_{c_+,0} = V_{c',0}$ and by the Picard-Lefschetz theorem,
$V_{c_+,1} = \tau_{V_{c',1}}(V_{c_-,1})$. Taking the vanishing cycles of $c''$
and moving them to $M$ by parallel transport shows that $HF(V_{c',1},V_{c_-,1})
= HF(V_{c',1},V_{c_+,1}) = HF(V_{c'',0},V_{c'',1})$. Applying Theorem
\ref{th:exact-sequence} in this context yields a long exact sequence
\[
\xymatrix{
 {HF(c'') \otimes HF(c')} \ar[r] &
 {HF(c_-)} \ar[d] \\
 & \ar@/^/[ul] {HF(c_+)}.
 }
\]
This can be thought of as a skein rule, with $HF(c'') \otimes HF(c')$ measuring
the change to $HF(c_-)$ as the path moves over $z \in S^\crit$ and becomes
$c_+$. For $S = D$, allowing oneself for a moment to believe naively that ``two terms in a
long exact sequence determine the third one'', one sees that all groups $HF(c)$
could be determined from knowing just finitely many of them, through a process
of successively breaking up paths into shorter pieces.

\section{Grading\label{sec:grading}}

\subsection{}
The assumption $2c_1(M,L) = 0$ in Theorem \ref{th:exact-sequence} is very close
to the condition needed to equip Floer cohomology groups with $\Z$-gradings, so
one might just as well take advantage of it. Doing that requires some
preliminaries, which we now go through.

Let $M$ be an arbitrary symplectic manifold, and $\JJ_M$ the space of all
compatible almost complex structures. There is a canonical unitary line bundle
$\Delta_M \rightarrow \JJ_M \times M$ whose fibre at $(J,x)$ is
$\Lambda^n(TM_x,J)^{\otimes 2}$. To any Lagrangian submanifold $L \subset M$
one can associate a section $\det^2(TL)$ over $\JJ_M \times L$ of the
associated circle bundle $S(\Delta_M)$; and for any symplectic automorphism
$\phi$ there is a canonical isomorphism $\det^2(D\phi): \Delta_M \rightarrow
\Delta_M$, covering the map $\phi_* \times \phi$ from $\JJ_M \times M$ to
itself. A Maslov map is a trivialization $\delta_M: \Delta_M \rightarrow \C$.
Suppose that we have chosen such a $\delta_M$. A grading of $L \subset M$ is a
lift
\[
\xymatrix{
 &&& {\R} \ar[d]^{\exp(2\pi i \cdot)} \\
 {\JJ_M \times L} \ar@/^1pc/[urrr]^-{\tL}
 \ar[rr]^-{\;\det^2(TL)} && {S(\Delta_M)} \ar[r]^-{\delta_M} &
 {S^1.\!\!}
}
\]
Any two gradings differ by a locally constant function $L \rightarrow \Z$; we
write $\tL[\sigma]$ for $\tL-\sigma$, $\sigma \in \Z$. Similarly, a grading of
$\phi$ is a map $\tphi: \JJ_M \times M \rightarrow \R$ such that $\exp(2\pi i\,
\tphi(J,x)) = \delta_M(\textstyle{\det^2}(D\phi)(w))/\delta_M(w)$ for any $w
\in S(\Delta_M)$ in the fibre over $(J,x)$. When dealing with
connected manifolds $M$
with boundary and maps $\phi$ trivial near $\partial M$, there is often a
preferred grading, characterized by being zero near $\JJ_M \times \partial M$.
Gradings of $L$ and $\phi$ induce a grading of $\phi(L)$:
\[
\tphi(\tL) \stackrel{\mathrm{def}}{=} (\tphi + \tL) \circ (\phi_*^{-1} \times
\phi^{-1}).
\]
If $L \subset M$ is a framed Lagrangian sphere admitting gradings, there is a
distinguished grading $\ttau_L$ of $\tau_L$, which is zero outside $\JJ_M
\times \{\text{neighbourhood of $L$}\}$ and satisfies
\begin{equation} \label{eq:self-shift}
\ttau_L(\tL) = \tL[1-n].
\end{equation}
Graded Lagrangian submanifolds and graded symplectic automorphisms are pairs
consisting of an $L$ resp.\ $\phi$ together with a choice of grading. For
brevity we denote such pairs by $\tL,\tphi$ only. The Floer cohomology of a
pair of graded Lagrangian submanifolds, whenever defined (e.g.\ when $M$ is
exact with contact type boundary, and the Lagrangian submanifolds are
themselves exact), has a canonical $\Z$-grading with the properties
\begin{align}
 \notag
 & HF^*(\tL_1[-\sigma],\tL_2) = HF^*(\tL_1,\tL_2[\sigma]) = HF^{*+\sigma}(\tL_1,\tL_2), \\
 \label{eq:graded-floer}
 & HF^*(\tphi(\tL_1),\tphi(\tL_2)) \iso HF^*(\tL_1,\tL_2), \\
 \notag
 & HF^*(\tL_2,\tL_1) \iso HF^{n-*}(\tL_1,\tL_2)^\vee.
\end{align}

\subsection{}
Now consider an exact Morse fibration $(E,\pi)$ over $S$, and the space
$\JJ_{E/S}$ of pairs $(j,J)$ such that $J$ is compatible relative to $j$. Let
$\Delta_{E/S} \rightarrow \JJ_{E/S} \times E$ be the line bundle with fibres
$\Lambda^{n+1}(TE_x,J)^{\otimes 2} \otimes (TS_{\pi(x)},j)^{\otimes -2}$. A
{\bf relative Maslov map} is a trivialization $\delta_{E/S}: \Delta_{E/S}
\rightarrow \C$. Once such a map has been chosen, there is a canonical induced
Maslov map $\delta_M$ on any regular fibre $M = E_z$. The vanishing cycles $V
\subset M$ associated to paths $c: [0;1] \rightarrow S$, $c(0) = z$ and $c(1)
\in S^\crit$, admit gradings; and the monodromies $\rho_l$ along closed loops
$l$ in $S \setminus S^\crit$, $l(0) = z$, even have canonical gradings, which
we denote by $\tilde{\rho}_l$. If $[l] \in H_1(S)$ vanishes, $\tilde{\rho}_l$
is zero near $\JJ_M \times \partial M$. This implies a graded version of the
Picard-Lefschetz theorem, saying that $\tilde{\rho}_l = \ttau_V$ for $c,l$ as
in Figure \ref{fig:doubling}.

\subsection{}
Any Lagrangian boundary condition $Q$ for $(E,\pi)$ comes with a canonical
section of $S(\Delta_{E/S})|\JJ_{E/S} \times Q$. We call $\delta_{E/S}$ adapted
to $Q$ if there is a diagram
\[
\xymatrix{
 {\JJ_{E/S} \times Q} \ar[d]_{\id \times \pi}
 \ar[rrr]^-{\text{canonical section}} &&&
 {S(\Delta_{E/S})} \ar[d]^{\delta_{E/S}} \\
 {\JJ_{E/S} \times \partial S} \ar@{-->}[rrr]^-{\lambda} &&&
 {S^1.\!\!}
}
\]
For compact $S$, the existence of an adapted relative Maslov map means that all
components of the space of $(j,J)$-holomorphic sections taking
$\partial S$ to $Q$ have the same expected
dimension $n\,\chi(S) + \deg(\lambda)$; we have already seen an instance of
this in Remark \ref{rem:technical}. There is a generalization of this notion to
the relative case, which ensures that the Floer groups in \eqref{eq:phirel}
have canonical gradings and that the invariant $\Phirel$ is of a specific
degree, not necessarily zero. Instead of explaining this in general we will
just look at an example, that of the long exact sequence.

Suppose then that we have an exact symplectic manifold $M$ with a Maslov map
$\delta_M$, graded Lagrangian submanifolds $\tL_1,\tL_2$, and a framed exact
Lagrangian sphere $L$ which admits gradings. The map $b: HF^*(\tL_1,\tL_2)
\rightarrow HF^*(\tL_1,\ttau_L(\tL_2))$ has degree zero, and the same holds for
the pair-of-pants product $a: HF^*(\tL,\tL_2) \otimes HF^*(\tL_1,\tL)
\rightarrow HF^*(\tL_1,\tL_2)$. The coproduct $HF^*(\tL_1,\ttau_L(\tL_2))
\rightarrow HF^*(\tL,\ttau_L(\tL_2)) \otimes HF^*(\tL_1,\tL)$ has degree $n$
because of Poincar{\'e} duality, but then by \eqref{eq:self-shift} and
\eqref{eq:graded-floer} the isomorphism
\[
HF^*(\tL,\ttau_L(\tL_2)) = HF^*(\ttau_L^{-1}(\tL),\tL_2) =
HF^{*+1-n}(\tL,\tL_2)
\]
has degree $1-n$. Hence the third map $c$ raises degrees by one, like in any
other cohomological long exact sequence.

\section{Mutation\label{sec:mutation}}

\subsection{}
As in our discussion of Floer theory, we take the ground field to be $\Z/2$;
all categories will be linear over it. Let ${\mathcal C}$ be a triangulated
category such that the spaces $\Hom_{\mathcal C}^*(X,Y) = \bigoplus_r
\Hom_{\mathcal C}(X,Y[r])$ are finite-dimensional for any $X,Y$. An exceptional
collection in ${\mathcal C}$ is a finite family of objects $(X^1,\dots,X^m)$
satisfying
\[
 \Hom_{\mathcal C}^*(X^i,X^k) =
 \begin{cases}
 \Z/2 \cdot \id_{X^i} & i = k, \\
 0 & i > k.
 \end{cases}
\]
Such a collection is called full if the $X^i$ generate ${\mathcal C}$ in the
triangulated sense. For $X,Y \in \Ob\,{\mathcal C}$ define $T_X(Y)$, up to
isomorphism, as the object fitting into an exact triangle
\[
T_X(Y)[-1] \rightarrow \Hom^*_{\mathcal C}(X,Y) \otimes X \xrightarrow{ev} Y
\rightarrow T_X(Y);
\]
here the tensor product is just a finite sum of shifted copies of $X$, and $ev$
the canonical evaluation map. If $(X^1,\dots,X^m)$ is a full exceptional
collection then so are $(Y^1,\dots,Y^m)$, $(Z^1,\dots,Z^m)$ where
\begin{align} \label{eq:first-mutation}
 Y^{i} &= \begin{cases}
 T_{X^{1}}(X^{i+1}) & i<m, \\
 X^{1} & i=m,
 \end{cases} \\
\intertext{respectively} \label{eq:second-mutation}
 Z^{i} &= \begin{cases}
 X^{i} & i<m-1, \\
 T_{X^{m-1}}(X^{m}) & i=m-1, \\
 X^{m-1} & i=m.
 \end{cases}
\end{align}

\subsection{}
There is a slightly different version of the same story for
$A_\infty$-categories. Call an $A_\infty$-category $\A$ {\bf directed} if it
has finitely many objects numbered $1,\dots,m$, say $\Ob\,\A =
\{X^1,\dots,X^m\}$, such that
\[
hom_\A(X^i,X^k) = \begin{cases}
 \text{finite-dimensional} & i<k, \\
 \Z/2 \cdot \id_{X^i} & i = k, \\
 0 & i>k.
\end{cases}
\]
Note that $\mu^d_\A = 0$ is necessarily zero for $d>\max\{m-1,2\}$. We recall the
definition of the (bounded) derived category $D^b(\A)$. The first step is to
embed $\A$ into a larger $A_\infty$-category $\Aoplus$ which has finite sums
and shifts. Thus, an object of $\Aoplus$ is a formal sum $\bigoplus_{e \in E}
X_e[\sigma_e]$ with $E$ a finite set, $X_e \in \Ob\,\A$, and $\sigma_e \in \Z$.
Next, a twisted complex in $\A$ is a pair $(C,\delta_C)$ consisting of $C \in
\Ob\,\Aoplus$ and $\delta_C \in hom^1_\Aoplus(C,C)$, such that the
``generalized Maurer-Cartan equation''
\begin{equation} \label{eq:maurer-cartan}
\sum_{d\geq 1} \mu_{\Aoplus}^d(\delta_C,\dots,\delta_C) = 0
\end{equation}
holds. Twisted complexes form an $A_\infty$-category $\Tw\,\A$ which again has
direct sums and shifts. It also contains a cone $\Cone(a)$ for any morphism $a$
such that $\mu^1_{\Tw\,\A}(a) = 0$, and the cohomological category $D^b(\A) =
H^0(\Tw\,\A)$ inherits a triangulated structure from this. The objects of $\A$,
seen as twisted complexes with zero differential, form a full exceptional
collection in $D^b(\A)$.

\begin{remark}
The definition of the derived category of a general $A_\infty$-category $\A$
uses pairs $(C,\delta_C)$ such that $\delta_C$ is strictly decreasing with
respect to some finite filtration of $C$, which has the effect of making the
sum \eqref{eq:maurer-cartan} finite. The fact that this is not necessary in our
case is just one of several technical simplifications which directedness brings
with it.
\end{remark}

Any $A_\infty$-functor $F: \A \rightarrow \B$ between directed
$A_\infty$-categories induces another one $\Tw\,F: \Tw\,\A \rightarrow \Tw\,\B$
taking cones to cones, and therefore an exact functor $D^b(F) = H^0(\Tw\, F):
D^b(\A) \rightarrow D^b(\B)$. Call $F$ a quasi-isomorphism if $H(F): H(\A)
\rightarrow H(\B)$ is an isomorphism; since the objects of $\A$ generate
$D^b(\A$), it follows that in this case $D^b(F)$ is an equivalence.

More interestingly, suppose that $\A$ is some directed $A_\infty$-category, and
$(Y^1,\dots,Y^m)$ an exceptional collection in $D^b(\A)$. One can then define a
new directed $A_\infty$-category $\B$, called the directed
$A_\infty$-subcategory of $\Tw\,\A$ generated by the $Y^i$, as follows. Objects
of $\B$ are the $Y^i$, and $hom_\B(Y^i,Y^k) = hom_{\Tw\,\A}(Y^i,Y^k)$ for
$i<k$, with the same composition maps between these groups as in $\Tw\,\A$. All
other morphism groups and compositions are as dictated by the directedness
condition. The embedding $\B \hookrightarrow \Tw\,\A$ can be extended to an
$A_\infty$-functor $\iota: \Tw\,\B \rightarrow \Tw\,\A$, which induces an exact
functor $H^0(\iota): D^b(\B) \rightarrow D^b(\A)$. For the same reason as
before, $H^0(\iota)$ is full and faithful; it is even an equivalence if
$(Y^1,\dots,Y^m)$ is a full exceptional collection.

\subsection{}
Given $X \in \Ob\,\Tw\,\A$ and a finite-dimensional complex $V$ of vector
spaces, one can form the tensor product $V \otimes X \in \Ob\,\Tw\,\A$, which
is a direct sum of shifted copies of $X$ with a differential combining those on
$V$ and $X$. Taking $V = (hom_{\Tw\,\A}(X,Y),\mu^1_\A)$ for some $Y \in
\Ob\,\Tw\,\A$, one has a canonical evaluation morphism $ev \in
hom_{\Tw\,\A}^0(hom_{\Tw\,\A}(X,Y) \otimes X, Y)$ with $\mu^1_{\Tw\,\A}(ev) =
0$. Let $T_X(Y) \in \Ob\,\Tw\,\A$ be the cone of $ev$. This is isomorphic in
$D^b(\A)$ to the object of the same name introduced above, but it is now unique
in a strict sense, not just up to isomorphism. Taking the exceptional
collection formed by the objects of $\A$ and applying \eqref{eq:first-mutation}
or \eqref{eq:second-mutation} yields another full exceptional collection, hence
a directed $A_\infty$-subcategory $\B$ of $\Tw\,\A$ with $D^b(\B) \iso
D^b(\A)$. This process can be repeated indefinitely and leads to the following
notion:

\begin{definition} \label{def:mutation}
Two directed $A_\infty$-categories with $m$ objects are mutations of each other
if they can be related by a sequence of the following moves and their inverses:
\begin{myitemize}
\item
$\A \rightsquigarrow \B$ if there is a quasi-isomorphism between them.

\item
It is allowed to shift each object by some degree, which means changing the
grading of each group $hom_\A(X^i,X^k)$ by $(\sigma_i-\sigma_k)$ for some
$\sigma_1,\dots,\sigma_m \in \Z$, while keeping the same composition maps.

\item
$\A \rightsquigarrow c\A$ where, if the objects of $c\A$ are denoted by
$\{Y^1,\dots,Y^m\}$, the nontrivial morphism spaces are
\[
\qquad hom_{c\A}(Y^i,Y^k) =
\begin{cases}
 hom_\A(X^{i+1},X^{k+1}) & \text{$i<k<m$}, \\
 hom_\A(X^{1},X^{i+1})^\vee[-1] & \text{$i<m,\, k=m$.}
\end{cases}
\]
Here ${}^\vee$ denotes the dual of a graded vector space. The compositions
$\textstyle{\mu^d_{c\A}: \prod_{\nu=1}^d hom_{c\A}(Y^{i_\nu},Y^{i_{\nu+1}})
\rightarrow hom_{c\A}(Y^{i_1},Y^{i_{d+1}})}$, $i_1 < \dots < i_{d+1}$, are
equal to those in $\A$ except when $i_{d+1} = m$, in which case one has
$\leftsc \mu^d_{c\A}(a^d,\dots,a^1),b \rightsc = \leftsc a^d,
\mu^d_\A(a^{d-1},\dots,a^1,b) \rightsc$ with $\leftsc \dots \rightsc$ the dual
pairing.

\item
$\A \rightsquigarrow r\A$ with $\Ob\,r\A = \{Z^1,\dots,Z^m\}$ and the following
nontrivial morphism spaces $hom_{r\A}(Z^i,Z^k)$:
\begin{equation} \label{eq:rmutation}
 \qquad \begin{cases}
 hom_\A(X^i,X^k) & \text{$i<k\leq m-2$,} \\
 hom_\A(X^i,X^{m-1}) & \text{$i \leq m-2$, $k = m$,} \\
 hom_\A(X^{m-1},X^m)^\vee[-1] & \text{$i = m-1$, $k = m$,} \\
 (hom_\A(X^{m-1},X^m) \otimes hom_\A(X^i,X^{m-1}))[1]
 \;\oplus \hspace{-6em} & \\
 \hspace{2em} \oplus \; hom_\A(X^i,X^m)
 & \text{$i \leq m-2$, $k = m-1$.}
 \end{cases}
\end{equation}
$\mu^1_{r\A}$ is given by $\mu^1_\A$ in the first two cases, its dual
$(\mu^1_\A)^\vee$ in the third, and in the final case by
\[
\begin{pmatrix}
 \mu^1_\A \otimes \id + \id \otimes \mu^1_\A & 0 \\
 \mu^2_\A & \mu^1_\A
\end{pmatrix}.
\]
Most of the higher order maps
\begin{equation} \label{eq:composition}
\mu^d_{r\A}: \textstyle{\prod_{\nu=1}^d hom_{r\A}(Z^{i_\nu},Z^{i_{\nu+1}})}
\rightarrow hom_{r\A}(Z^{i_1},Z^{i_{d+1}})
\end{equation}
for $i_1 < \dots < i_{d+1}$ are taken from those of $\A$ in a straightforward
way, but there are two exceptions. One is when $i_{d+1} = m-1$, in which case
\[
 \mu^d_{r\A} =
 \begin{pmatrix}
 \id \otimes \mu^d_\A & 0 \\
 \mu_\A^{d+1} & \mu^d_\A
 \end{pmatrix}
\]
with respect to the obvious splittings on both sides of \eqref{eq:composition}.
The second exceptional case is when $i_d = m-1$, $i_{d+1} = m$: then
$\mu^2_{r\A}$ is
\[
\quad\,
\begin{array}{c}%
 hom_\A(X^{m-1},X^m)^\vee
 \otimes hom_\A(X^{m-1},X^m) \otimes
 hom_\A(X^{i_1},X^{m-1})
 \\
 \oplus\;\, (hom_\A(X^{m-1},X^m)^\vee \otimes
 hom_\A(X^{i_1},X^m))[-1]
 \\
 \Big\downarrow 
 \\
 hom_\A(X^{i_1},X^{m-1}), \vspace{0.5em}
 \\
 \mu^2_{r\A}(a^3 \otimes a^2 \otimes a^1, b^2 \otimes b^1) =
 \leftsc a^3, a^2 \rightsc a^1,
 \end{array}
\]
while the compositions of order $d \geq 3$ vanish.
\end{myitemize}
\end{definition}

Actually, while $r\A$ is precisely the directed $A_\infty$-subcategory of
$\Tw\,\A$ generated by the collection \eqref{eq:second-mutation}, $c\A$ is only
canonically quasi-isomorphic to that generated by \eqref{eq:first-mutation}.
But it is still true that two directed $A_\infty$-categories which are
mutations of each other have equivalent derived categories.

\section{Fukaya categories}

\subsection{}
Throughout this section $M$ is a fixed exact symplectic manifold, with a Maslov
map $\delta_M$. A graded Lagrangian configuration in $M$ is a family
$\widetilde{\Gamma} = (\tL_1,\dots,\tL_m)$ of graded, exact, framed Lagrangian
spheres. Hurwitz equivalence for graded configurations is defined as in the
ungraded case, with the following adaptations: for the isotopy invariance one
wants to take a grading $\tphi$ which is zero near $\JJ_M \times \partial M$,
so that $\tphi$ is isotopic to the identity in the group of graded symplectic
automorphisms; the moves $c\widetilde{\Gamma},r\widetilde{\Gamma}$ use the
canonical gradings $\ttau_L$ of Dehn twists; and there is an additional shift
move,
\begin{itemize}
\item
$\widetilde{\Gamma} \rightsquigarrow (\tL_1[\sigma_1],\dots,\tL_m[\sigma_m])$
for any $\sigma_1,\dots,\sigma_m \in \Z$.
\end{itemize}
This is something of an anticlimax, since it cancels out the extra information
contained in the grading; but in fact, the whole notion of graded configuration
has been introduced only for notational convenience.

We will associate to $\widetilde{\Gamma}$ a {\bf directed Fukaya
$A_\infty$-category} $\Lag(\widetilde{\Gamma})$, unique up to
quasi-isomorphism. Suppose first that the configuration is in general position,
meaning that any two $L_i$ are transverse and there are no triple
intersections. Objects of $\A = \Lag(\widetilde{\Gamma})$ are the graded
Lagrangian submanifolds $\tL_i$, in the given order, and
\[
hom_\A(\tL_i,\tL_k) =
\begin{cases}
 CF^*(\tL_i,\tL_k) = (\Z/2)^{L_i \cap L_k} & i<k, \\
 \Z/2 \cdot \id_{\tL_i} & i=k, \\
 0 & i>k.
\end{cases}
\]
Roughly speaking, $\mu^1_\A$ is the Floer boundary map, $\mu^2_\A$ the
pair-of-pants product, and $\mu^3_\A,\mu^4_\A,\dots$ Fukaya's generalizations
of that product. Each $\mu^d_\A$ depends on the choice of a family $\J^{d+1}$
of almost complex structures on $M$, and these choices have to obey certain
consistency conditions, which we will now outline (this is joint work of
Lazzarini and the author).

In a first step one takes, for each $1 \leq i_1 < i_2 \leq m$, a generic family
of almost complex structures $\J^2(i_1,i_2) : [0;1] \rightarrow \JJ_M$. As is
well known, this causes all solutions of Floer's equations
\begin{equation} \label{eq:floer}
 u: \R \times [0;1] \rightarrow M, \qquad
 \left\{
  \begin{split}
  & \J^2(i_1,i_2,t) \circ du(s,t) = du(s,t) \circ j, \\
  & u(\R \times \{1\}) \subset L_{i_1}, \; u(\R \times \{0\}) \subset L_{i_2}, \\
  & \textstyle{\int u^*\o} < \infty,
  \end{split}
 \right.
\end{equation}
where $j$ is the standard complex structure on $\R \times [0;1]$, to be
regular. From the one-dimensional solution spaces (zero-dimensional after
dividing by translation) one builds $\mu^1_\A: CF^*(\tL_{i_1},\tL_{i_2})
\rightarrow CF^{*+1}(\tL_{i_1},\tL_{i_2})$.

Next take $S = D$, three cyclically ordered marked points $\zeta_\nu^3 \in
\partial S$, $1 \leq \nu \leq 3$, local coordinates $\psi^3_\nu: D^+
\rightarrow S$ around them, and set $S^* = S \setminus
\{\zeta_1^3,\zeta_2^3,\zeta_3^3\}$, all as in the definition of the invariants
$\Phirel$. Choose for each $1 \leq i_1 < i_2 < i_3 \leq m$ a generic family
$\J^3(i_1,i_2,i_3): S^* \rightarrow \JJ_M$, such that for $s \ll 0$ and $t \in
[0;1]$,
\begin{equation} \label{eq:tubular}
\J^3(i_1,i_2,i_3,\psi^3_\nu(e^{\pi(s + i t)})) =
\begin{cases}
 \J^2(i_{\nu},i_{\nu+1},t) & \nu = 1,2, \\
 \J^2(i_1,i_3,1-t) & \nu = 3.
\end{cases}
\end{equation}
Denote by $I^3_\nu$, $1 \leq \nu \leq 3$, the connected components of $\partial
S^*$, ordered cyclically such that $I_1^3$ lies between $\zeta_3^3$ and
$\zeta_1^3$. The equation defining $\mu^2_\A$ is
\[
 u: S^* \rightarrow M, \qquad
 \left\{
  \begin{split}
  & \J^3(i_1,i_2,i_3,z) \circ du(z) = du(z) \circ j, \\
  & u(I_\nu^3) \subset L_{i_\nu} \qquad \text{for $\nu = 1,2,3$,} \\
  & \textstyle{\int u^*\o} < \infty.
  \end{split}
 \right.
\]
Condition \eqref{eq:tubular} causes this to agree with a suitable equation
\eqref{eq:floer} in the tubular coordinates $\psi^3_\nu(e^{\pi(s+it)})$ on each
end of $S^*$.

The additional ingredient in the definition of the products of order $d > 2$
are ``moduli parameters'' as the complex structure of the domain changes. Let
$\CC^{d+1} \subset (\partial D)^{d+1}$ be the configuration space of $d+1$
distinct, numbered and cyclically ordered points on $\partial D$. The moduli
space $\RR^{d+1}$ and the universal disc bundle $\SS^{d+1}$ over it are defined
as
\[
 \SS^{d+1} = \CC^{d+1} \times_{\Aut(D)} D
 \longrightarrow \RR^{d+1} = \CC^{d+1}/\Aut(D),
\]
where $\Aut(D) \iso PSL(2,\R)$ is the holomorphic automorphism group. Each
fibre $\SS^{d+1}_r$ carries a canonical complex structure, and there are
canonical sections $\zeta^{d+1}_\nu: \RR^{d+1} \rightarrow \partial\SS^{d+1}$,
$1 \leq \nu \leq d+1$, such that the points $\zeta^{d+1}_\nu(r) \in \partial
\SS^{d+1}_r$ are distinct and cyclically ordered for each $r$. Write
\[
\SS^{d+1,*} = \SS^{d+1} \setminus (\textstyle{\bigcup_\nu}
\im\,\zeta_\nu^{d+1}), \qquad \SS^{d+1,*}_r = \SS^{d+1,*} \cap \SS^{d+1}_r.
\]
For each $1 \leq i_1 < \dots i_{d+1} \leq m$ one has to choose a family
$\J^{d+1}(i_1,\dots,i_{d+1}): \SS^{d+1,*} \rightarrow \JJ_M$, subject to two
kinds of conditions.
\begin{myitemize}
\item
Compatibility with $\J^2$. This requires a preliminary choice of maps
$\psi^{d+1}_\nu: \RR^{d+1} \times D^+ \rightarrow \SS^{d+1}$, $1 \leq \nu \leq
d+1$, such that $\psi^{d+1}_\nu(r,\cdot)$ provides local coordinates around
$\zeta^{d+1}_\nu(r)$ for each $r \in \RR^{d+1}$. Then the conditions are
similar to \eqref{eq:tubular}, requiring
$\J^{d+1}(i_1,\dots,i_{d+1},\psi^{d+1}_\nu(r,z))$ to be determined by the
previously chosen $\J^2$ for small $|z|$.

\item
Compatibility with $\J^{e+1}$ for $2 \leq e < d$. For this it is necessary to
consider the Deligne-Mumford compactification of $\RR^{d+1}$. Each stratum at
infinity is a product of lower order spaces $\RR^{e+1}$, and for a point $r \in
\RR^{d+1}$ sufficiently close to one such stratum, the fibre $\SS^{d+1,*}_r$ is
built by gluing together fibres of $\SS^{e+1,*}$ for the various occurring $e$.
The precise condition on $\J^{d+1}(i_1,\dots,i_{d+1})|\SS^{d+1,*}_r$ is too
complicated to be written down here, but informally it says that this should be
built up from the $\J^{e+1}$ in a corresponding way.
\end{myitemize}
Let $I^{d+1}_{r,\nu}$, $1 \leq \nu \leq d+1$, be the connected components of
$\partial\SS^{d+1,*}_r$, ordered cyclically so that $I^{d+1}_{r,1}$ lies
between $\zeta^{d+1}_{d+1}(r)$ and $\zeta^{d+1}_1(r)$. The consistency
conditions leave enough freedom to make solutions of the equation
\[
r \in \RR^{d+1},\, u: \SS^{d+1,*}_r \rightarrow M, \quad \left\{
\begin{split}
 & \J^{d+1}(i_1,\dots,i_{d+1},z) \circ du(z) = du(z) \circ j, \\
 & u(I^{d+1}_{r,\nu}) \subset L_{i_\nu} \quad \text{for $\nu = 1,\dots,d+1$}, \\
 & \textstyle{\int} u^*\o < \infty.
\end{split}
\right.
\]
regular for generic $\J^{d+1}(i_1,\dots,i_{d+1})$. Counting such solutions
defines
\[
\mu^d_\A: CF^*(\tL_{i_d},\tL_{i_{d+1}}) \otimes \cdots \otimes
CF^*(\tL_{i_1},\tL_{i_2}) \longrightarrow CF^{*+2-d}(\tL_{i_1},\tL_{i_{d+1}}).
\]
The dependence of $\Lag(\widetilde{\Gamma})$ on the choice of almost complex
structure can be analyzed using a one-parameter family argument. Since only
finitely many moduli spaces are involved, the $A_\infty$-structure is subject
to a finite number of changes in the family. At each of these exceptional times
one can produce a quasi-isomorphism relating the old $A_\infty$-structure with
the new one.

\begin{remark}
As a technical point, note that directedness allows us to bypass some
problems which plague Fukaya's original setup, having to do with the chain
complexes underlying $HF(L,L)$ and unit elements in them.
\end{remark}

The next step is isotopy invariance which, as always in Floer theory, is also
used to extend the definition to configurations which are not in general
position.

\begin{prop} \label{th:isotopy-invariance}
Let $\widetilde{\Gamma} = (\tL_1,\dots,\tL_m)$ be a graded Lagrangian
configuration in general position. Take $l \in \{1,\dots,m\}$, a symplectic
automorphism $\phi$ isotopic to the identity, and a grading $\tphi$ which is
zero near $\JJ_M \times \partial M$, such that $\widetilde{\Xi} =
(\tL_1,\dots,\tL_{l-1}, \tphi(\tL_l),\dots,\tphi(\tL_m))$ is again in general
position. Then there is a quasi-isomorphism $F: \Lag(\widetilde{\Gamma})
\rightarrow \Lag(\widetilde{\Xi})$.
\end{prop}

We will spend a moment discussing the structure of the proof, since it is a
good example of arguments involving directed Fukaya categories. Recall that an
$A_\infty$-functor $F: \A \rightarrow \B$ consists of a map $F: \Ob\,\A
\rightarrow \Ob\,\B$, chain maps $F^1: hom_\A(X,Y) \rightarrow hom_\B(FX,FY)$
for $X,Y \in \Ob\,\A$, and multilinear ``higher order terms'' $F^d$, $d \geq
2$. In the present case, the map on objects is the obvious one, and the
nontrivial chain maps
\[
F^1: CF^*(\tL_i,\tL_k) \longrightarrow CF^*(\tL_i,\tphi(\tL_k)), \quad i<l
\text{ and } k \geq l,
\]
are those underlying the continuation homomorphisms, so that $F$ is
automatically a quasi-isomorphism. The main effort goes into defining higher
order terms which satisfy the equations for an $A_\infty$-functor.

\subsection{}
We will now describe the relation between Hurwitz moves of $\widetilde{\Gamma}$
and mutations of $D^b\Lag(\widetilde{\Gamma})$. Proposition
\ref{th:isotopy-invariance} says that isotopies of $\widetilde{\Gamma}$ result
in a quasi-isomorphism. Shifting the gradings $\tL_i$ obviously corresponds to
the first mutation in Definition \ref{def:mutation}.

\begin{lemma} \label{th:cyclic}
$\Lag(c\widetilde{\Gamma})$ is quasi-isomorphic to $c\Lag(\widetilde{\Gamma})$.
\end{lemma}

The proof relies on the $\Z/(d+1)$-action on $\RR^{d+1}$ given by a cyclic
shuffle of the marked points (this symmetry had not been used in the definition
of directed Fukaya categories).

\begin{theorem} \label{th:main}
$r\Lag(\widetilde{\Gamma})$ is quasi-isomorphic to $\Lag(r\widetilde{\Gamma})$.
\end{theorem}

To see why this is plausible, set $\A = \Lag(\widetilde{\Gamma})$, and let
$(Z^1,\dots,Z^m)$ be the objects of $r\A$. The cohomology
$H(\hom_{r\A}(Z^i,Z^k),\mu^1_{r\A})$, $i < k$, is
\[
 \qquad \begin{cases}
 HF^*(\tL_i,\tL_k) & \text{$i<k \leq m-2$,} \\
 HF^*(\tL_i,\tL_{m-1}) & \text{$i \leq m-2$, $k = m$,} \\
 HF^*(\ttau_{L_{m-1}}(\tL_m),\tL_{m-1}) & \text{$i = m-1$, $k = m$,} \\
 H(\Cone(\mu^2_\A: CF^*(\tL_{m-1},\tL_m) \otimes CF^*(\tL_i,\tL_{m-1})
 \;\rightarrow \hspace{-6em} & \\
 \hspace{2em} \rightarrow \; CF^*(\tL_i,\tL_m)))
 & \text{$i \leq m-2$, $k = m-1$.}
 \end{cases}
\]
This is just \eqref{eq:rmutation} except that in writing down the third case we
have used Poincar{\'e} duality and \eqref{eq:self-shift}. As we saw when
discussing Theorem \ref{th:exact-sequence}, the cone in the last line is
isomorphic to $HF^*(\tL_i,\ttau_{L_{m-1}}(\tL_m))$. Therefore all cohomology
groups are in fact isomorphic to the corresponding ones in
$\Lag(r\widetilde{\Gamma})$. The remainder of the proof, as in Proposition
\ref{th:isotopy-invariance}, consists in extending this to a full-fledged
$A_\infty$-functor.

By the general theory of mutations, what we have shown implies that if two
graded Lagrangian configurations in $M$ are Hurwitz equivalent, their directed
Fukaya categories have equivalent derived categories. Combining this with
Picard-Lefschetz theory yields the following consequence:

\begin{corollary} \label{th:invariant}
Let $(E,\pi)$ be an exact Morse fibration over $D$, with a relative Maslov map
$\delta_{E/D}$. Make an admissible choice of paths, and let $\Gamma$ be the
corresponding distinguished basis of vanishing cycles in a fibre. Choose any
gradings $\widetilde{\Gamma}$ and form $D^b\Lag(\widetilde{\Gamma})$. This is
independent of all choices up to equivalence, and hence is an invariant of
$(E,\pi)$ and $\delta_{E/D}$.
\end{corollary}

\subsection{}
There is a computational aspect which Corollary \ref{th:invariant} fails to
convey, and which we will explain by giving an example. Let $M$ be an exact
symplectic four-manifold with $2c_1(M) = 0$, and $(L_1,L_2)$ two Lagrangian
spheres in $M$. $L_2' = \tau_{L_1}^2(L_2)$ and $L_2$ are always isotopic as
smooth submanifolds. There are two cases, $L_1 = L_2$ and $L_1 \cap L_2 =
\emptyset$, in which $L_2'$ is also Lagrangian isotopic to $L_2$ for obvious
reasons, but in general this is false:

\begin{prop} \label{th:knotted}
Suppose that $L_1, L_2$ intersect transversally, with $|L_1 \cap L_2| \geq 3$,
and that the local intersection numbers at all points are the same. Then $L_2'$
is not Lagrangian isotopic to $L_2$.
\end{prop}

The proof goes as follows. Choose a Maslov map $\delta_M$ and gradings
$\tL_1,\tL_2$. The directed Fukaya category $\A = \Lag(\widetilde{\Gamma})$ of
the configuration $\widetilde{\Gamma} = (\tL_1,\tL_1,\tL_2,\tL_2)$ is
determined up to quasi-isomorphism by the graded vector space $R =
HF^*(\tL_1,\tL_2)$ together with the degree two maps $q_1,q_2 \in
\mathrm{End}(R)$ given by the pair-of-pants product with the unique nontrivial
element in $HF^2(\tL_1,\tL_1)$ resp.\ $HF^2(\tL_2,\tL_2)$. These satisfy $q_1
\circ q_2 = q_2 \circ q_1$ and $q_1^2 = q_2^2 = 0$. Lemma \ref{th:cyclic} and
Theorem \ref{th:main} give an explicit sequence of mutations which transforms
$\A$ into the directed Fukaya category $\B$ associated to the Hurwitz
equivalent configuration
\[
ccrc^{-1}rc^{-1}\widetilde{\Gamma} = (\ttau_{L_1}^2(\tL_2),\tL_1,\tL_1,\tL_2).
\]
In particular this determines $HF^*(\ttau_{L_1}^2(\tL_2),\tL_2)$, since that is
a morphism space in $H(\B)$. We will not write down the actual computation; the
outcome is the total cohomology of the complex
\[
 \Z/2 \oplus \Z/2[-2] \xrightarrow{(\id,q_2)}
 \mathrm{End}(R) \xrightarrow{\psi}
 \mathrm{End}(R)[2],
\]
where the second arrow is $\psi(x) = q_1 \circ x - x \circ q_1$. Since
$\psi^2(x) = 2 q_1 \circ x \circ q_1 = 0$, linear algebra tells us that the
dimension of $\coker(\psi)$ is $\geq (\dim R)^2/2$. With the assumption $\dim R
\geq 3$ this implies that $HF(L_2',L_2)$ is bigger than $HF(L_2,L_2)$, which
completes the argument.

\providecommand{\bysame}{\leavevmode\hbox to3em{\hrulefill}\thinspace}

\end{document}